\documentclass[12pt,a4paper,draft]{amsart}
\usepackage[english]{babel} 
\usepackage{a4wide,amssymb,xypic} 
\newtheorem{thm}{Theorem}[section] 
\newtheorem{corol}[thm]{Corollary} 
\newtheorem{lemma}[thm]{Lemma} 
\newtheorem{prop}[thm]{Proposition} 
\newtheorem{defin}[thm]{Definition} 
 
\theoremstyle{remark} 
\newtheorem{remark}[thm]{Remark} 
\newtheorem{example}[thm]{Example}
 
\newcommand{\Oc}{\mathcal O} 
\newcommand{\Ec}{\mathcal E}

\newcommand{\PP}{{\mathbb P}}
\newcommand{\F}{\mathcal{F}}

\newcommand{\R}{\mathbb R} 
\newcommand{\C}{\mathbb C} 
\newcommand{\Q}{\mathbb Q}

\def\Sym{\operatorname{Sym}} 
\def\rk{\operatorname{rk}} 
\def\coker{\operatorname{coker}} 
 
\def\dim{\operatorname{dim}} 
 
\newcommand\grass{\mbox{Gr}} 
\newcommand\hgrass{{\mathfrak{Gr}}} 
\newcommand{\cO}{{\mathcal O}}

\newcommand{\fE}{{\mathfrak E}}
\newcommand{\fF}{{\mathfrak F}}
\newcommand{\fG}{{\mathfrak G}}
\newcommand{\fQ}{{\mathfrak Q}}
\newcommand{\OPQ}[1]{{\mathcal O_{\PP Q_{#1}}}}
\begin{document} 
\title[Semistability vs. nefness]{SEMISTABILITY vs.  NEFNESS FOR (HIGGS) VECTOR BUNDLES} 
\date{7 June 2004, revised 18 March 2005}
\thanks{Research partly supported by the Spanish {\sc dges} through the 
research project BFM2003-0097,  by ``Junta de Castilla y Le\'on'' 
through the research project SA118/04, and by the Italian National Research Project ``Geometria delle variet\`a algebriche.''  Both authors are members 
of the research group {\sc vbac} (Vector Bundles on Algebraic
Curves).}
\maketitle  \thispagestyle{empty}
\begin{center}{\sc U. Bruzzo} \\ 
Scuola 
Internazionale Superiore di Studi Avanzati,\\ Via Beirut 2-4, 34013 
Trieste, Italia\\ E-mail {\tt bruzzo@sissa.it} \\[10pt]
{\sc D. Hern\'andez Ruip\'erez} \\
Departamento de Matem\'aticas, Universidad de Salamanca, 
\\ Plaza de la Merced 1-4, 37008 Salamanca, Espa\~na\\ E-mail {\tt ruiperez@usal.es}\end{center} 

\bigskip 
 
\begin{quote} 
\begin{center}{\sc Abstract} \end{center}\smallskip \small  Generalizing a result of Miyaoka, we prove that the semistability  of a vector bundle $E$ on a smooth projective curve over
a field of characteristic zero is equivalent to the nefness of any  of certain divisorial classes 
$\theta_s$, $\lambda_s$ in the Grassmannians  $\grass_s(E)$ of locally-free quotients of $E$ and in the projective bundles  $\PP Q_s$, respectively (here 
$0<s<\rk E$ and $Q_s$ is the universal quotient bundle on $\grass_s(E)$). The  result is   extended to Higgs bundles. In that case a necessary and sufficient condition for semistability is that  all classes $\lambda_s$ are nef. We also extend this result
to higher-dimensional complex projective varieties by showing that the nefness
of the classes $\lambda_s$ is equivalent to the semistability of the bundle $E$
together with the vanishing of the characteristic class $\Delta(E)=c_2(E)-\frac{r-1}{2r}c_1(E)^2$.
\end{quote} 
 
\newpage

\section{Introduction} 
Let   $E$ be a vector bundle on a smooth projective curve $C$ over a field 
of characteristic 0. According to Miyaoka \cite{Mi}, the semistability of 
$E$ is equivalent to the fact that 
  a suitable divisor $\lambda$ in the projective bundle $\PP E$ 
is numerically effective (nef) \cite{Mi}. The  
numerical class $\lambda$ is defined as 
$$\lambda = [c_1(\Oc_{\PP E}(1))] - \frac1r(\deg E)\,F$$ 
where $r=\rk E$ and $F$ is class of the fibre of the projection $\PP E\to C$. 
 
A mild generalization of Miyaoka's criterion is the following. We consider the Grassmannians $\grass_s(E)$ 
of rank $s$ locally-free quotients of $E$  and  the 
projectivized universal quotient bundles $\PP Q_s$ on them. Moreover we introduce the 
classes 
$$\theta_s=c_1(Q_s)-\frac sr\deg(E)\,F_s$$ 
in $\grass_s(E)$ and 
$$\lambda_s=c_1(\Oc_{\PP Q_s}(1)) - \mu(E) \, F_s$$ 
in $\PP Q_s$ (in both cases $F_s$ denotes the class of the fibre of  
the projection 
onto $C$). We prove: 
 
\begin{thm} If $E$ is semistable, all classes $\theta_s$ and  
$\lambda_s$ are nef. 
Conversely, if one of these classes is nef (for just one value of $s$),  $E$ is semistable. 
\label{primothm}\end{thm} 
 
The interest in this result is that it suggests a generalization of Miyaoka's criterion to the case of  Higgs bundles $(E,\phi)$ on a smooth projective curve $C$. 
(Throughout this paper we shall use the notation $\fE$ for a Higgs bundle
$(E,\phi)$.)
We introduce  schemes $\hgrass_s(\fE) $ parametrizing  locally-free 
rank $s$ Higgs quotients of $\fE $ (i.e., locally-free quotients of 
$E$ whose kernels are $\phi$-invariant), and the universal quotient Higgs 
bundles 
$\fQ_s=(Q_s,\Phi_s)$ on them. We 
define classes $\theta_{s,\fE}$  in $\hgrass_s(\fE)$ and $\lambda_{s,\fE}$ in $\PP Q_s$ as in the previous 
case. 
 
\begin{thm} If $\fE $ is semistable, all classes $\theta_{s,\fE}$  
and $\lambda_{s,\fE}$ are nef. 
Conversely, if \emph{all}   classes $\lambda_{s,\fE}$ are nef, then  
$\fE $ is semistable. 
\label{secthm}\end{thm} 
 
It is indeed not enough that one of the classes  $\lambda_{s,\fE}$ is nef 
to ensure the semistability of $\fE $, as we show in Section \ref{s3}. 

In Example \ref{exMiyaoka} we apply this criterion to the Higgs bundle
whose semistability implies the Miyaoka inequality for the Chern classes
of a projective surface with ample canonical bundle. 

The last generalization we give is about Higgs bundles on
a complex projective manifold $X$ of any dimension. (We need
to restrict to the complex case as  we use transcendental techniques.) In this case we define the classes 
$$\lambda_{s,\fE}=c_1(\OPQ{s}(1))-\frac1r \pi^\ast_s c_1(E)$$
where $\pi_s\colon \hgrass_s(\fE)  \to X$ is the projection. Let
$\Delta(E)$ be the characteristic class
$$\Delta(E) =  c_2(E)-\frac{r-1}{2r}c_1(E)^2 = \frac{1}{2r}c_2(E\otimes E^\ast)\, .$$
We prove the following result. 

\begin{thm} \label{thirdthm} Let $\fE$ be a Higgs bundle
on a complex projective manifold. The following conditions are equivalent.
\begin{itemize} \item[i)] All classes $\lambda_{s,\fE}$ are nef,
for $0<s<r$. \item[ii)]
$\fE $ is semistable and $\Delta(E)=0$.
\end{itemize}
\end{thm}

Note that   the classes $\lambda_{s,\fE}$ do not
depend on the choice of a polarization in $X$, so when they are all
nef the Higgs bundle $\fE $ is semistable with respect to all polarizations. Conversely,
if $\fE $ is semistable with respect to a given polarization, and  $\Delta(E)=0$,
then $\fE $ is semistable with respect to {\it all} polarizations. Our results
also imply that if all classes $\lambda_{s,\fE}$ are nef, then $\fE $ is semistable after restriction to {\it any} smooth projective curve in $X$.

Theorem \ref{thirdthm} applies
also to ordinary bundles in the following form:

\begin{thm} Let $E$ be a vector bundle on a complex projective manifold.
The following conditions are equivalent.
\begin{itemize} \item[i)] The class $\lambda_1$ is nef. \item[ii)]
$E$ is semistable and $\Delta(E)=0$.
\end{itemize}
\end{thm}

This result was already contained in \cite{BB} as a special case
and is proved by repeating verbatim the proof of Theorem \ref{thirdthm}
(see Section \ref{4}).

\medskip 
 
\noindent{\bf Acknowledgments.} We thank M.S.~Narasimhan for useful  
discussions and the referee for helping us to improve the presentation.
This paper was partly 
written while the first author was visiting the Tata Institute 
for Fundamental Research in Mumbai, 
to which thanks are due for hospitality and support.

\medskip 
 
\section{Semistability vs. nefness for vector bundles\label{s2}} 
All varieties we shall 
consider will be over an algebraically closed field of characteristic  
0.  Let $X$ be a 
smooth projective variety, $E$ a holomorphic vector bundle on it, and 
denote by $\PP E$ the projectivization of $E$, defined as 
$$\PP E = \mbox{\bf Proj}(\mathcal S(E)),$$ where $\mathcal 
S(E)$ is the symmetric algebra of the sheaf of sections of $E$. We recall 
that a bundle $E$ is said to be 
\emph{ample} if the hyperplane line bundle $\Oc_{\PP E}(1)$ on $\PP E$ 
is ample  \cite{Hart1,Hart2}. A weaker notion is 
that of \emph{numerical effectiveness:} a bundle $E$ is said to be 
numerically effective (nef) if  the class 
$c_1(\Oc_{\PP E}(1))$ is numerically effective. If both $E$ and $E^\ast$
are numerically effective, then $E$ is said to be \emph{numerically flat.}
The following result has been proved in \cite{DPS}.

\begin{prop}\label{vanishchern} 
All Chern classes of a numerically flat bundle vanish.
\qed\end{prop}
 
We recall that given a smooth projective variety $X$ with a choice of a 
polarization $H$, a torsion-free coherent sheaf $\Ec$ on $X$  is said to be 
\emph{semistable} (in Mumford-Takemoto's sense) if for every proper 
coherent subsheaf $\F$ of $\Ec$ one has 
$$\mu(\F)\le\mu(\Ec),$$ where the slope $\mu(\Ec)$ of a 
torsion-free coherent sheaf is defined as 
$$\mu(\Ec)=\frac{\deg\Ec}{\rk\Ec},\qquad \deg\Ec=c_1(\Ec)\cdot H^{n-1}$$ 
if $n=\dim X$. If the inequality always holds strictly, the bundle $E$ is 
said to be \emph{stable.} 
 
Following Miyaoka \cite{Mi}, we introduce    some notation 
and state the theorem relating the semistability of a vector bundle 
to the nefness of a suitable divisorial class. 
\begin{defin} For every smooth projective  variety $X$, one denotes: 
\begin{enumerate} 
\item $N^1(X)=\displaystyle\frac{\mbox{\rm Pic}(X)}{\mbox{num.  
eq.}}\otimes \R$; \smallskip 
\item $NA(X)\subset N^1(X)$ the ample cone of $X$ (the cone generated 
by the classes of ample divisors), and $\overline{NA}(X)$ its 
closure (the set of classes of nef divisors of $X$); \smallskip 
\item $N_1(X)=\displaystyle\frac{A_1(X)}{\mbox{num. eq.}}\otimes \R$;  
\smallskip 
\item $NE(X)\subset N_1(X)$ the real cone generated by the effective 1-cycles. 
\end{enumerate} 
\end{defin} 
\begin{thm}{\rm (\cite{Mi}, Theorem 3.1.)} Let  $\pi\colon E\to C$ be a rank $r$  
vector bundle on a smooth 
projective curve $C$,  and define the class in $N^1(\PP E)$ 
$$\lambda = [c_1(\Oc_{\PP E}(1)]-\mu(E)\,F$$ 
where $F$ is the class in $N^1(\PP E)$ of the fibre of the projection  
$\PP E\to C$. Then the following conditions 
are equivalent:
\begin{enumerate} 
\item $E$ is semistable;
\item $\lambda$ is nef;
\item $\overline{NA}(\PP E)=\R_+\,\lambda +\R_+\,F$; 
\item  $\overline{NE}(\PP E)=\R_+\,\lambda^{r-1} +\R_+\,\lambda^{r-2}\cdot F$; 
\item every effective divisor in $\PP E$ is nef.\hfill $\square$ 
\end{enumerate}\label{Mithm}\end{thm} 
\noindent (Here $\R_+$ is the set of nonnegative real numbers.) The class 
$r\lambda$ is the relative anti-canonical class of $\PP E\to C$, and 
one has $\lambda^r=0$. 
\begin{proof} For the readers' convenience, and following \cite{Mi},
we include a sketch of a proof of this theorem.

(i) $\Rightarrow$ (ii). If $\lambda$ is not nef there is an irreducible curve
$C'\subset \PP E$ such that $[C']\cdot\lambda<0$. After a suitable base change 
$f\colon C'' \to C$ we may assume that $C'$ is a union of sections $C_j$ of the
bundle $\PP (f^\ast E)$, and $[C_j]\cdot \lambda'<0$ for all $j$, where
$\lambda'$ is the class $\lambda$ for the bundle $\pi''\colon \PP(f^\ast E) \to C''$. 
There are surjections $f^\ast E \to \pi''_\ast \Oc_{C_j}(1)$, and
$\deg  (\Oc_{C_j}(1)) = [C_j]\cdot\lambda' + \mu(f^\ast E) < \mu(f^\ast E)$.
But this contradicts the semistability of $E$ (note that since the morphism
$f$ is separable and finite, the bundle $E$ is semistable if and only if $f^\ast E$ is,
cf.~\cite[Prop. 3.2]{Mi} and our Lemma \ref{l:pullback}). 

(ii)  $\Rightarrow$ (iv). If $\Gamma$ is a class in $\overline{NE}(\PP E)$, one has
$\Gamma=a\lambda^{r-1}+b\lambda^{r-2}\cdot F$ with $a\ge 0$. Since $\lambda$ is nef,
one has $ b = \Gamma\cdot\lambda \ge 0$. 

(iii) and (iv) are easily shown to be equivalent.

(iii) and (iv) $\Rightarrow$ (v). Let $D=a\lambda+b\, F$ 
be an effective divisor class. All 1-cycles $D\cdot (\lambda +\varepsilon F)^{r-2}$
lie in $\overline{NE}(\PP E)$ for every positive real number $\varepsilon$,
and so do their limits $D\cdot \lambda^{r-2}$. Then $a,b\ge0$ by (iv) and
$D$ is nef by (iii). 

(v) $\Rightarrow$ (i). Let $F$ be a destabilizing subbundle of $E$ and let
$\alpha\in\Q$ be such that $\mu(F)>\alpha>\mu(E)$. For $N$ big enough
the space
\begin{eqnarray*} H^0(\Sym\!^NF(-N\alpha p)) &\subset&
H^0(C,\Sym\!^NE(-N\alpha p)) \\ &\simeq& H^0(\PP E,\Oc_{\PP E}(N) \otimes \pi^\ast \Oc_C (-N\alpha p)) 
\end{eqnarray*}
(where $p$ is a point in $C$) is nonempty; therefore, the class $N(\lambda+(\mu(E)-\alpha)F))$ is effective
but not nef.
\end{proof}
 
We describe now the first generalization of this result. 
Given a vector bundle $E$ on an algebraic variety $X$, 
we shall denote by $\grass_s(E)$ the Grassmann variety  
of rank $s$ 
locally-free quotients of $E$, with $ 0 < s < r=\rk E$. We have a morphism 
$p_s\colon\grass_s(E)\to X$ that makes $\grass_s(E)$ a bundle of Grassmannians. On every 
variety  $\grass_s(E)$ a universal quotient bundle $Q_s$ is defined, 
in a such a way that for any morphism $f\colon Y\to X$ and any rank $s$ 
locally-free quotient $F$ of $f^\ast E$, there is a morphism 
$\psi_F\colon Y\to\grass_s(E)$ over $X$ (that is, $f=p_s\circ\psi_F$) 
such that $F=\psi_F^\ast Q_s$. 
 
Let $\theta_s$ be the class in $N^1(\grass_s(E))$ 
$$\theta_s=[c_1(Q_s)]-\frac sr\deg(E)\,F_s,$$ 
where $F_s$ is the class of the fibre of the projection $\pi_s\colon \grass_s(E)\to X$. 
\begin{thm}\label{firstcon} If $E$ is a semistable vector bundle on a smooth projective curve $C$, the class $\theta_s$ 
is nef for every $s$, $0<s<r=\rk E$. 
\end{thm} 
\begin{proof} This result can be proved 
according to the lines of the proof of the implication (i) $\Rightarrow$ 
(ii) in Theorem \ref{Mithm}.
Alternatively, one can use the Pl\"ucker embedding $\varpi\colon \grass_s(E)  
\to \PP(\Lambda^s E)$ to reduce to the case of a  
projective bundle, as we now show. The  morphism $\varpi$  embeds into 
a commutative diagram 
$$\xymatrix{ \PP Q_s \ar[r] \ar[d] & \PP(\Lambda^s E) \ar[d]^{\pi_\Lambda} 
\\ 
\grass_s(E) \ar[ur]^\varpi \ar[r]^{p_s} & C } $$ 
and one has an isomorphism 
\begin{equation}\label{tuttovabene} \varpi^\ast \Oc_{ \PP(\Lambda^s  
E)}(1)\simeq \det Q_s. 
\end{equation} 
The induced morphisms $$\varpi^\ast\colon  
N^1(\PP(\Lambda^s E))\to N^1(\grass_s(E))\,,\qquad\varpi_\ast\colon 
N_1(\grass_s(E))\to N_1(\PP(\Lambda^s E))$$ are isomorphisms, as one 
 easily shows  by using some Schubert calculus. 
 
If $E$ is  
semistable, so is $\Lambda^s E$ 
for every $s$, $0<s<r$ \cite{Maru}, 
so that the class 
$$\lambda_{\Lambda^s E}=[c_1(\Oc_{ \PP(\Lambda^s  
E)}(1)]-\mu(\Lambda^s E) \,[\pi_\Lambda^\ast(x)]$$ 
is nef. By restricting to the image of the Pl\"ucker embedding one  
obtains that $\theta_s$ is nef.  
\end{proof} 
 
The converse to Theorem \ref{firstcon} is as follows. 
\begin{thm} If for some $s$ (with $0<s<r=\rk E$) the class $\theta_s$  
is nef, then 
$E$ is semistable. 
\label{tnis}\end{thm} 
\begin{proof} If $\theta_s\in N^1(\grass_s(E))$ is nef,   
the class $\lambda_{\Lambda^s E}=(\varpi^\ast)^{-1}(\theta_s)$ 
is nef as well, since for any curve $\Gamma\subset \PP(\Lambda^s E)$
one has $ \lambda_{\Lambda^s E}\cdot [\Gamma] = \theta_s \cdot [\Gamma']$
with $\varpi(\Gamma')=\Gamma\cap \varpi(\grass_s(E))$.
By Miyaoka's result
the bundle $\Lambda^s E$ is semistable. It is then an easy task to  
prove that $E$ is semistable as well. 
\end{proof} 
 
\begin{remark}  
These constructions provide an alternative algebraic proof of the fact  
that, given a semistable bundle $E$ 
on a smooth projective variety $X$, its exterior powers $\Lambda^sE$  
are semistable as well. By the 
Metha-Ramanathan theorem (cf.~e.g.~\cite{Mi}), it is enough to  
consider the case when $X$ 
is a curve. Then by the proof of Theorem \ref {firstcon} which follows 
\cite{Mi} we know that 
the classes $\theta_s$ are nef, so that the classes  
$\lambda_{\Lambda^s E}$ are nef, whence 
$\Lambda^sE$  is semistable.\qed\end{remark} 
 
We consider now another construction. Again, $E$ is a rank $r$ vector bundle 
on a smooth projective curve $C$, $\grass_s(E)$ is the Grassmannian bundle of 
its rank $s$ quotients, and $Q_s$ the universal quotient bundle on  
$\grass_s(E)$. We define the 
class in $N^1(\PP Q_s)$ 
$$\lambda_s=[c_1(\Oc_{\PP Q_s}(1))] - \mu(E) \, F_s$$ 
where $F_s$ is the class of the fibre of the composition 
$ \PP Q_s \to \grass_s(E) \to C$. 
 
\begin{thm}\label{div} If $E$ is semistable, the class $\lambda_s$ 
is nef for every $s$, $0<s<r=\rk E$.\hfill$\square$ 
\end{thm} 
To prove this result we need a Lemma. 
\begin{lemma}\label{lemmino} Let $G$ be a    semistable vector bundle  
on a smooth 
projective curve $C$, and let $C'$ be an irreducible curve in $\PP  
G$. Denote by $\xi$ the 
class of $\Oc_{\PP G}(1)$ in $N^1(\PP G)$. Then, 
$$[C']\cdot\xi \ge \mu(G)\,p_\ast[C']$$ 
where $p\colon \PP G\to C$ is the projection. 
\end{lemma} 
\begin{proof} 
By Theorem \ref{Mithm} we have 
$$[C']=a\lambda^{r-1}+b\lambda^{r-2}\cdot F$$ 
(here $\lambda = [c_1(\Oc_{\PP G}(1))-\mu(G) F)]$ and $r=\rk G$) 
with $a,\,b\ge 0$, so that 
$$[C']\cdot\xi =  
(a\lambda^{r-1}+b\lambda^{r-2}F)(\lambda+\mu(G)F)=a\mu(G)+b\ge  
a\mu(G).$$ 
Moreover, one has $a=p_\ast[C']$. 
\end{proof} 
\noindent \emph{Proof of Theorem \ref{div}.} If for some $s$ the  
class $\lambda_s$ is not nef 
there is an irreducible curve $C'\subset \PP Q_s$ which surjects onto  
$C$ and is such that 
$C'\cdot\lambda_s<0$. Let $h\colon C''\to C$ be a finite morphism and  
consider the 
commutative diagram whose squares are cartesian 
$$\xymatrix{ 
\PP Q'_s \ar[r]^{h''} \ar[d]_{\pi'_s} & \PP Q_s \ar[d]^{\pi_s} \\ 
\grass_s(h^\ast E) \ar[r]^{h'} \ar[d]_{p'_s} &\grass_s(E)  \ar[d]^{p_s} \\ 
C'' \ar[r]^h & C 
}$$ 
We may choose the pair $(C'',h)$ in such a way that the fibre  
product $\tilde C = C''\times_C C'$ 
(a curve in $\PP Q'_s$) is a union of curves $C_j$ which project onto  
$C''$ with degree one (and meet 
the fibre $F'_s$ at just 
one point). 
One has $[C_j]\cdot \lambda'_s<0$. 
Let $\Gamma_j$ be the projection of $C_j$ onto $\grass_s(h^\ast E)$,  
denote by $Q_j$ 
the restriction of $Q'_s$ to it, and let $E_j=({h'}^\ast\circ  
p_s^\ast E)_{\vert \Gamma_j}$. We have an epimorphism 
$E_j \to Q_j\to 0.$ 
The composition $h\circ p'_s$ restricted to $\Gamma_j$ (call it  
$h_j$) is a finite morphism (actually, 
an isomorphism), so 
that $E_j$ is semistable. 
Now with the help of Lemma \ref{lemmino} we have: 
\begin{eqnarray*}\mu(Q_j) &\le & \displaystyle\frac{[C_j]\cdot  
\xi_j}{{\pi'_s}_\ast [C_j]} = [C_j]\cdot \xi'_s \\[3pt] 
&=&  [C_j]\cdot \left(\lambda'_s+\mu (h^\ast E)\,F'_s\right) < \mu  
(E_j)\end{eqnarray*} 
but this contradicts the semistability of $E_j$.\hfill$\square$ 
 
\begin{corol}  \label{amplequot} If $E$ is a semistable bundle of  
positive (resp.~nonnegative) 
degree on a smooth projective curve $C$, 
then all universal quotient bundles $Q_s$ are ample (resp.~numerically effective). 
\end{corol} 
\begin{proof}  This result may   be proved   by mimicking Gieseker's  
proof for $s=1$, cf.~\cite{Gies}. 
  \end{proof}

\begin{thm}  If for some $s$ (with $0<s<r=\rk E$) the class  
$\lambda_s$ is nef, then 
$E$ is semistable. 
\label{lnis}\end{thm} 
\begin{proof}  By direct computation one sees that  
$\pi_{s\ast}(\lambda_s)^s=\theta_s$. So, if $\Gamma$ is 
a curve in $\grass_s(E)$, one has  
$[\Gamma]\cdot\theta_s=(\lambda_s)^s\cdot \pi_s^\ast [\Gamma]\ge 0$ 
(since $\lambda_s$ is nef) so that $\theta_s$ is nef, whence $E$ is semistable. 
\end{proof}

  \medskip 
 
\section{Semistability vs. nefness for Higgs bundles\label{s3}} 
 
We want to investigate if the semistability of a Higgs bundle 
can be encoded in the nefness of some suitable classes. In particular, we prove 
Theorem \ref{secthm}.  
 
\subsection{Grassmannians of Higgs quotients} We recall the basic definitions 
about Higgs bundles (cf.~\cite{Simp1}, \cite{Simp2}). 
 
\begin{defin} Let $X$ be a projective variety. A Higgs sheaf $\fE$ on $X$ is a 
coherent sheaf $E$ on 
$X$ endowed with a morphism $\phi\colon E\to E\otimes \Omega_X$ of 
$\Oc_X$-modules such that 
$\phi\wedge\phi=0$, where $\Omega_X$ is the cotangent sheaf to $X$. 
A Higgs subsheaf $F$ of a Higgs sheaf $\fE=(E,\phi)$ is a subsheaf of $E$
such that $\phi(F)\subset F\otimes\Omega_X$.
A  Higgs bundle is a Higgs sheaf 
$\fE $ such that $E$ is a locally-free $\Oc_X$-module. 
\end{defin}

\begin{defin} Let $X$ be a smooth projective variety equipped with a 
polarization. A Higgs sheaf $\fE=(E,\phi) $ is semistable (resp.~stable) if it 
is torsion-free, and $\mu(F)\le \mu(E)$ (resp. $\mu(F)< \mu(E)$) for 
every proper nontrivial Higgs subsheaf $F$ of $\fE$.\end{defin}

In the sequel we shall need the following Lemma. It generalizes 
a well-known fact about semistable vector bundles that we 
have already used in this paper \cite{Gies,Mi}. 
 
\begin{lemma}  Let $f\colon Y\to X$ be a finite separable morphism of 
smooth projective curves, $\fE$ a Higgs  bundle on $X$ and $f^\ast\fE$ the pullback Higgs bundle on $Y$. Then $\fE $ is semistable 
if and only if $f^\ast \fE$ is semistable. 
\label{l:pullback} 
\end{lemma} 
 
\begin{proof} The ``if'' part is straightforward. To prove that $f^\ast \fE$ is semistable
  on $Y$ when $\fE $ is semistable on $X$ we can  assume that $f$ is a
Galois covering.
 Let $0\to F\to f^\ast E$ be a maximal destabilizing  Higgs subbundle.  For every element $\sigma$ in the Galois group of $f$, $\sigma^\ast F$ is also a maximal destabilizing Higgs
subbundle of
  $\sigma^\ast E=E$ so that $\sigma^\ast F=F$ by the uniqueness of $F$. It
follows that $F=f^\ast E'$ for a certain
   subbundle $0\to E'\to E$. Since $E'$ destabilizes $E$,  we have only to
prove that
  it is a Higgs subbundle of $\fE $, i.e., that the composed morphism
    $E' \xrightarrow{\phi_{\vert E'}} E\otimes \Omega_X\to
(E/E')\otimes \Omega_X$
vanishes. Since $f$ is faithfully flat, we may as well prove that the
induced morphism
$F \xrightarrow{f^\ast\phi} f^\ast E\otimes f^\ast \Omega_X\to (f^\ast
E/F)\otimes f^\ast\Omega_X$ vanishes. But this follows from  the diagram
$$
\xymatrix{ & &0\ar[d] &0\ar[d]& \\ & F\ar[r]^{f^\ast \phi\quad
}\ar[rd]^{\phi_Y} & f^\ast E\otimes f^\ast \Omega_X \ar[d] \ar[r] &
(f^\ast E/F)\otimes f^\ast\Omega_X \ar[d] & \\ 0\ar[r] & F\otimes
\Omega_Y \ar[r] & f^\ast E\otimes \Omega_Y \ar[r] & (f^\ast E/F)\otimes
\Omega_Y \ar[r]  &0 }
$$
since $\phi_Y\colon F\to  f^\ast E\otimes \Omega_Y$ takes values in
$F\otimes \Omega_Y$.
\end{proof}

Given a Higgs sheaf $\fE $, we may construct the closed subschemes 
$\hgrass_s(\fE)\subset \grass_s(E)$ parametrizing  
the rank $s$ locally-free Higgs  quotients, i.e. 
locally-free quotients of $E$ such that the corresponding kernels are 
$\phi$-invariant. This can be done as follows: let us consider the  
universal exact sequence 
\begin{equation} 
0\to S_{r-s}\xrightarrow{\psi} p_s^\ast E\xrightarrow{\eta} Q_s\to 0 
\label{eq:univ} 
\end{equation} 
of sheaves on the Grassmannian $\grass_s(E)$ that defines the  
universal quotient bundle $Q_s$. Then  $\hgrass_s(\fE)$ is the closed  
subvariety of $\grass_s(E)$ where the composed morphism 
$$ 
(\eta\otimes1)\circ p_s^\ast(\phi) \circ \psi\colon S_{r-s}\to Q_s\otimes 
 p_s^\ast\Omega_X 
$$ 
vanishes (the equations for $\hgrass_s(\fE)$ inside 
$\grass_s(E)$ are written in subsection \ref{eqns} for the special 
case $s=1)$). We denote by $\pi_{s,\fE}$ the projections 
$\hgrass_s(\fE)\to X$. 
The restriction of \eqref{eq:univ} to the scheme $\hgrass_s(\fE)$ 
 gives a new  universal exact sequence 
$$ 
0\to S_{r-s}\xrightarrow{\psi} (\pi_{s,\fE})^\ast E\xrightarrow{\eta} 
 Q_{s,\fE}\to 0 
$$ 
and $Q_{s,\fE}$ is a rank $s$ universal Higgs  quotient vector bundle.  
This means that  for every morphism $f\colon Y\to X$ and every  rank $s$ Higgs 
quotient $F$ of $f^\ast E$ there is a morphism $\psi_F\colon Y\to 
\hgrass_s(\fE)$ such that $f=\pi_{s,\fE}\circ \psi_F$ and $F=\psi_F^\ast Q_{s,\fE}$. Note 
that the kernel $S_{r-s}$ 
of the morphism $(\pi_{s,\fE})^\ast E \to Q_{s,\fE}$ is $\phi$-invariant.

For every $s$ we define the classes 
$$\theta_{s,\fE}\in N^1(\hgrass_s(\fE)),\qquad \lambda_{s,\fE}\in  
N^1(\PP Q_{s,\fE})$$ 
as in the previous Section. We have: 
\begin{thm} If $\fE $ is a semistable Higgs bundle on a smooth 
projective curve $C$, all   classes $\theta_{s,\fE}$ and $\lambda_{s,\fE}$ are nef. 
\label{sinhb}\end{thm} 
\begin{proof} A possible proof of the nefness of $\theta_{s,\fE}$ runs  
as in the proof 
of Theorem \ref{firstcon} which follows \cite{Mi}. Analogously, the proof of the nefness of  
$\lambda_{s,\fE}$ runs as in the 
proof of Theorem \ref{div}. 
\end{proof} 

\begin{remark} The proof of Theorem \ref{div} adapted to the case
of Higgs bundles shows that if
$\lambda_{s,\fE}$ is not nef, after a base change $\fE $  is destabilized
by a rank $s$ locally free quotient.
\end{remark}

\begin{corol}  If $\fE $ is a semistable Higgs bundle on a smooth 
projective curve $C$ of positive (resp.~nonnegative) degree,    then 
for all $s$ the universal 
quotient bundle $Q_s$ is ample (resp.~numerically effective). 
\end{corol} 
\begin{proof} One again adapts the proof by Gieseker in \cite{Gies}, this 
time using Lemma \ref{l:pullback}. 
\end{proof}

\subsection{Equations of the scheme of rank-one Higgs quotients\label{eqns}} 
For some time we concentrate on the scheme $\hgrass_1(\fE)$ 
of rank one Higgs quotients of a 
rank $r$ Higgs vector bundle  $\fE $  on a $n$-dimensional smooth  
variety $X$.   
We have a closed immersion $j \colon 
\hgrass_1(\fE)\hookrightarrow \PP E$, and the universal Higgs quotient 
is $Q_1=\Oc_{\hgrass_1(\fE)}(1)=j^\ast  \Oc_{\PP E}(1)$. We denote by 
$\pi_1\colon \hgrass_1(\fE)\to X$ the projection. 
 
We denote by $\xi_E=c_1(\Oc_{\PP E}(1))$ the hyperplane class in 
$\PP E$ and write 
$\xi_{\fE }=c_1(\Oc_{\hgrass_1(\fE)}(1))$ $=j^\ast(\xi_E)$. 
We also define 
$$ 
\lambda_{\fE }=\left[\xi_{\fE }-\frac1r\pi^{\ast}_1c_1  E\right]\in 
N^1(\hgrass_1(\fE))\,. 
$$ 
 
One can write local equations for $\hgrass_1(\fE)$ by using the  Euler   sequence 
$$ 0 \to\Omega_{\PP (E)/X}\otimes \Oc_{\PP E}(1) \xrightarrow{\ \psi\ } 
\pi^\ast E \xrightarrow{\ \eta\ } \Oc_{\PP E}(1) \to 0\,, 
$$ 
because this is the form that \eqref{eq:univ} takes in this case. We then know that 
$\hgrass_1(\fE)$    is the closed set where the  composition of morphisms 
$$ (\eta\otimes 1)\circ\pi^\ast\phi\circ\psi \colon\Omega_{\PP (E)/X}\otimes 
\Oc_{\PP E}(1)\to \Oc_{\PP E}(1)\otimes \pi^\ast\Omega^1_X 
$$  vanishes. Given a local 
  basis of sections $(e_1,\dots,e_r)$ of $E$, which can be taken as local 
vertical homogeneous coordinates for $\PP E$,  the Higgs field is  
represented by a 
matrix 
$(\phi_{\alpha\beta})$ of 1-forms by letting 
$\phi(e_\beta)=\sum_\alpha\phi_{\alpha\beta}e_\alpha$. The  
homogeneous equations for 
$\hgrass_1(\fE)$ are 
\begin{equation} 
\sum_\gamma  
e_\gamma(\phi_{\gamma\beta}e_\alpha-\phi_{\gamma\alpha}e_\beta)=0\,,\mbox{  
for 
every $1\le \alpha<\beta\le r$.} 
\label{e:piphi} 
\end{equation} 
So $\hgrass_1(\fE)$ is locally the intersection of $n\binom r2$ 
hyperquadrics in $\PP E$. Let us study this locus  in the case   when the Higgs bundle $\fE $  is \emph{nilpotent},  
i.e., there is a 
decomposition 
$$ E=E_1\oplus\dots\oplus E_m 
$$   as a direct sum of subbundles, and $\phi(E_i)\subseteq 
E_{i+1}\otimes\Omega_X$ for $1\le i< m$,  $\phi(E_m)=0$. 
 
The induced morphism 
$\phi\otimes 1\colon E\otimes\Omega_X^\ast\to E$ 
yields a Higgs quotient sheaf  $\fQ=(Q,0)$ of $\fE $, where
$Q=\coker(\phi\otimes1)$,
and there is also a Higgs quotient  bundle  $\bar{\fE}=(E/E_m,\bar\phi)$, where  
$\bar\phi$ is the Higgs field induced by $\phi$. 
There are  closed immersions 
$$ 
\hgrass_1(\fQ)=\PP(Q)\hookrightarrow \hgrass_1(\fE)\,, 
\qquad 
\hgrass_1(\bar{\fE})\hookrightarrow \hgrass_1(\fE)\,. 
$$ 
The homogeneous equations for $\PP(Q)$ are given locally by the images of a basis of \break $(\phi\otimes 1)(E\otimes\Omega_X^\ast)$, that is, by 
\begin{equation}
\sum_\beta \phi_{\alpha\beta} e_\beta=0\, \qquad \text{for every $e_\alpha\in 
E_i$, $ i<m$. }
\label{e:q}
\end{equation} 
\begin{prop}  The scheme of rank one Higgs quotients is the closed  
subscheme of $\PP(E)$ given by 
$$ 
\hgrass_1(\fE)= \hgrass_1(\bar{\fE})\cup\hgrass_1(\fQ)\cup Z\,. 
$$ 
where $Z\subseteq \hgrass_1(\bar{\fE})\cap\hgrass_1(\fQ)$ is a union  
of embedded components. 
\label{p:nilpotentpiphi} 
\end{prop} 
\begin{proof}  We can proceed locally. Let us write $n_i=\dim (E_1\oplus\dots \oplus E_i)$ and take for every
subbundle $E_i$ a local basis  
$\{e_\gamma\}$ of sections with $n_{i-1}<\gamma\le n_i$. If we consider the subset of the  equations  \eqref{e:piphi}  where 
$e_\beta\in E_m$ and $e_\alpha\in E_i$ for $i<m$, we get  
$$ 0=e_\beta \cdot  \sum_{n_i<\gamma\le n_{i+1}} e_\gamma \phi_{\gamma\alpha}\,, 
$$ for every pair $(\alpha,\beta)$ as above. These equations describe the locus  
$Y\cup Y'$, where 
\begin{align*} Y' & \equiv \{ e_\beta=0 \,\vert\, \text{for every $e_\beta\in E_m$}\} \equiv
\PP(E/E_m)\\ Y &\equiv \{ 
\sum_{n_i<\gamma\le n_{i+1}} e_\gamma \phi_{\gamma\alpha}=0 \,\vert\, \text{for every $e_\alpha\in 
E_i$, $ i<m$ }\} \equiv 
\PP(Q)\,,
\end{align*} where the last equality is due to \eqref{e:q}. 
The remaining equations are 
$$ 0=e_\alpha\cdot\big(  \sum_{n_j<\gamma\le n_{j+1}} e_\gamma \phi_{\gamma\beta} \big) - 
e_\beta\cdot \big ( 
\sum_{n_i<\gamma\le n_{i+1}} e_\gamma \phi_{\gamma\alpha} \big) 
$$ for $e_\alpha\in E_i$, $e_\beta\in E_j$, $i\le j<m$ (and $\alpha<\beta$ if $i=j$) and 
define hyperquadrics containing $Y$.  
 
If $m=2$, the only possibility is $i=j=1=m-1$ so that the $e_\gamma$'s in the  
equations above belong to $E_m$, thus proving that the corresponding hyperquadric contains 
also $Y'$. Moreover, $\bar\phi=0$ in this case, so that $Y'=\hgrass_1(\bar{\fE})$ and we 
conclude that  
$$ 
\hgrass_1(\fE)=\hgrass_1(\bar{\fE})\cup \hgrass_1(\fQ)=\PP(E/E_2) \cup \PP(Q) 
$$ when $m=2$.
 
Assume now that $m\ge 3$. 
Since the points in $Y\equiv \PP(Q)$ satisfy the equations \eqref{e:piphi}, we have that
$$
\hgrass_1(\fE)=(\hgrass_1(\fE)\cap Y') \cup \PP(Q) 
$$
Moreover, the equations for $\hgrass_1(\fE)\cap Y'$ are the equations \eqref{e:piphi} for $j\le m-1$, which are easily shown to be the equations for $\hgrass_1(\bar{\fE})$. Then $\hgrass_1(\fE)\cap Y' \equiv \hgrass_1(\bar{\fE})$ and $\hgrass_1(\fE)=\hgrass_1(\bar{\fE}) \cup \PP(Q)$ up to a union of embedded components.
\end{proof} 

\subsection{Equations in the case of curves} Let $\fE $ be a nilpotent 
Higgs bundle on   a smooth projective curve $C$ and 
denote by $\widetilde\hgrass_1(\fE)$  the union of all components of $\hgrass_1(\fE)$ not 
contained in a fibre of $\PP(E)\to C$. Similar meaning will have the  
expressions 
 $\widetilde\hgrass_1(\fQ)$ or  $\widetilde\hgrass_1({\fE})$. The symbol
 $\bar\lambda_{1,\fE}$ will denote the restriction of $\lambda_{1,\fE}$ to 
 $\widetilde\hgrass_1(\bar{\fE})$.
\begin{prop} 
The class of $\widetilde\hgrass_1(\fE)$ in the Chow ring of $\PP(E)$ is 
\begin{align*} 
\lbrack \widetilde\hgrass_1(\fE)\rbrack &=[\widetilde\hgrass_1(\fQ)] +  
j_\ast[\widetilde\hgrass_1(\bar{\fE})]\\ 
&=  \xi^{r-r(\phi)} - [\deg(\phi(E))+r(\phi)(2-2 g(C))+\deg(T(Q))] 
\xi^{r-r(\phi)-1}\cdot F \\ 
  &\phantom{=\ }+ j_\ast[\widetilde\hgrass_1(\bar{\fE})] 
\end{align*} 
where $r(\phi)=\rk (\phi(E))$, $j\colon \PP(E/E_m)\hookrightarrow  
\PP(E)$ is the natural immersion and $T(Q)$ is the torsion subsheaf of $Q$.
\label{p:class} 
\end{prop} 
\begin{proof} 
We start by computing the class of $\PP(G)$ where 
$$ 
0\to N\to E\to G\to 0 
$$ is a quotient rank $q$ bundle. One has 
$$ 
[\PP(G)]=a \xi ^{r-q}+b\xi^{r-q-1}\cdot F 
$$ 
where $\xi$ is the relative hyperplane class of $\PP(E)$, $F$ 
is the class of a fibre of $\pi$ and $a$, $b$ are integer numbers.  
Since $1=\xi_G^{q-1}\cdot 
F_G=\xi^{q-1}\cdot F\cdot [\PP(G)]$ we obtain $a=1$. Moreover  
$\xi^r=\pi^\ast(c_1(E))\cdot 
\xi^{r-1}=r\mu(E) F\cdot \xi^{r-1}=r\mu(E)$ and similarly  
$\xi_G^q=q\mu(G)$, and we get 
$b=q\mu(G)-r\mu(E)$, that is 
\begin{equation} 
\begin{aligned} 
\lbrack \PP(G)\rbrack &= \xi^{r-q}+(q\mu(G)-r\mu(E)) \xi^{r-q-1}\cdot F \\ 
& 
= \xi^{r-q}-(\deg (E)-\deg (G)) \xi^{r-q-1}\cdot F= \xi^{r-q}-\deg  
(N) \xi^{r-q-1}\cdot F\,. 
\end{aligned} 
\label{e:quotclass} 
\end{equation} 
When $G$ has torsion $T(G)$,  we actually have $\widetilde\hgrass_1(\fQ)=\PP(G/T(G))$. Since 
$G/T(G)$ is a quotient vector bundle, we can compute as above to  get 
\begin{equation} 
\lbrack \widetilde\hgrass_1(\fQ)\rbrack=\lbrack \PP(G/T(G))\rbrack =  
\xi^{r-q}-(\deg (N)+\deg(T(G))) 
\xi^{r-q-1}\cdot F\,. 
\label{e:quotclass2} 
\end{equation} 
This formula implies our claim. 
\end{proof} 

\subsection{Unstable Higgs   bundles $\fE $ such that 
$\tilde\lambda_{1,\fE }$ is nef.} 
Let $\fE $ be a rank three nilpotent Higgs bundle on a smooth  
projective curve $C$, having 
the form  $E=L_1\oplus L_2\oplus L_3$ where 
each $L_i$ is a line bundle and 
$\phi(L_1)\subseteq L_2\otimes\Omega_C$, $\phi(L_2)\subseteq  
L_3\otimes\Omega_C$, 
$\phi(L_3)=0$. 
 
Let us write $\alpha_i=c_1(L_i)$.  If we impose that the Higgs  
subbundles $L_3$ and $L_2\oplus L_3$ do not 
destabilize $\fE $ we obtain the inequalities 
\begin{align} 
\alpha_1+\alpha_2-2\alpha_3 &\ge 0 \quad\text{and}\label{e:31}\\ 
2\alpha_1-\alpha_2-\alpha_3 & \ge 0 \label{e:32} 
\end{align} One can prove that these inequalities are actually sufficient for 
$\fE $ to be semistable. 
 
Now we  study when the restrictions of $\tilde\lambda_{\fE }$ to  
the  components   of $\widetilde\hgrass_1(\fE)$ 
are nef. There are  two components, 
$$ 
  \widetilde{\hgrass_1}(\fE) =  \widetilde{\hgrass_1}(\bar{\fE})\cup \widetilde{\hgrass_1}(\fQ) 
$$ 
with $ \widetilde{\hgrass_1}(\bar{\fE})\simeq \widetilde{\hgrass_1}(\fQ) \simeq C$.  
Here $\bar{\fE}$ is the Higgs bundle given by $E/L_3$ with the
induced Higgs morphism.
Let us write 
$\tilde\lambda_{1}$ and $\tilde\lambda_{2}$ for the restrictions of 
$\tilde\lambda_{1,\fE }$ to each of the components of $\widetilde\hgrass_1(\fE)$.
 
Since $E/L_3\simeq L_1\oplus L_2$, by the nilpotent rank  
two case we have 
\begin{align*} 
\lbrack \widetilde{\hgrass_1}(\bar{\fE})\rbrack & =2 (\xi -\alpha_2 F)\cdot \lbrack 
\PP(L_1\oplus L_2)\rbrack \\ & =2 (\xi -\alpha_2 F)(\xi -\alpha_3 F)) \\ &=2 
(\xi ^2-(\alpha_2+\alpha_3) \xi \cdot F)\,, 
\end{align*} so that 
$$ 
\tilde\lambda_{1}= 2(\xi^2-(\alpha_2+\alpha_3) 
\xi \cdot F)(\xi -\frac13(\alpha_1+\alpha_2+\alpha_3) 
F)=\frac23\,(2\alpha_1-\alpha_2-\alpha_3)\,. 
$$ 
 
On the other hand, since 
$$ Q=L_1\oplus (L_2/(\phi\otimes 1)(L_1\otimes\Omega_C^\ast)) \oplus 
(L_3/(\phi\otimes 1)(L_2\otimes\Omega_C^\ast)) 
$$ 
by modding   the torsion out we obtain $Q/T(Q)=L_1$ and 
$$ 
\lbrack \widetilde{\hgrass_1}(\fQ) \rbrack = \lbrack \PP(L_1)\rbrack= 
\xi^2-(\alpha_2+\alpha_3)\xi \cdot F\,, 
$$ by Eq \eqref{e:quotclass2}. Then 
\begin{align*} 
\tilde\lambda_{2} &= 
(\xi^2-(\alpha_2+\alpha_3)\xi\cdot F)(\xi  
-\frac13(\alpha_1+\alpha_2+\alpha_3)F) \\ & = 
\frac13\,(2\alpha_1-\alpha_2-\alpha_3)\,. 
\end{align*} So $\tilde\lambda_{1}$ and $\tilde\lambda_{2}$  
are nef if and only if 
inequality \eqref{e:32} holds, and if that inequality holds, 
$\tilde\lambda_{1,\fE }$ is nef.

Let $C$ be a smooth projective curve of genus 2, and $K=x+y$ a canonical 
divisor. Let us consider the line bundles 
$$ L_1=\Oc_C(K+x)\,,\quad      L_2=\Oc_C(x)\,,\quad   L_3=\Oc_C(3x) 
$$ Since $L_2\otimes L_1^{-1}\otimes \Omega_C=\Oc_C$, there exists a nonzero 
morphism $\phi_{21}\colon L_1\to L_2\otimes\Omega_C$. Moreover, $L_3\otimes 
L_2^{-1}\otimes \Omega_C=\Oc_C(2x+K)$ so that there is also a nonzero morphism 
$\phi_{32}\colon L_2\to L_3\otimes\Omega_C$. We can then define a  
nilpotent Higgs 
field $\phi\colon E\to E\otimes\Omega_C$ on 
$$ E=L_1\oplus L_2\oplus L_3 
$$ as being equal to $\phi_{21}$ on $L_1$, to $\phi_{32}$ on $L_2$ and zero on 
$L_3$. Now, 
$$ 2\alpha_1-\alpha_2-\alpha_3=2\,,$$ 
that is, the inequality  \eqref{e:32} is true, so that 
$\tilde\lambda_{1,\fE }$ is nef. The restriction of the class  
$\lambda_{1,\fE }$ 
to a component of $\hgrass_1(\fE)$ lying in a fibre of $\PP E\to C$ 
coincides with the restriction of the class $\xi_E=[c_1(\Oc_{\PP  
E}(1))]$, hence 
is nef, and the     class $\lambda_{1,\fE }$ itself is nef. 
However, 
$$ 
\alpha_1+\alpha_2 -2 \alpha_3 = -2 \,, 
$$ so that the inequality  \eqref{e:31} does not hold, and   $\fE $ is 
not semistable. 
 
It is interesting to check in this example what happens  with  the
class $\theta_{2,\fE}$ 
in $\hgrass_2(\fE)$. If we again remove the components embedded in  
fibres of $\PP E\to C$, 
we obtain a subscheme $\widetilde{\hgrass_2}(\fE)\simeq \PP  
L_3\simeq C$, and the 
class $\tilde\theta_{2,\fE}$ (the restriction of $\theta_{2,\fE}$ in   to 
  $\widetilde{\hgrass_2}(\fE)$) is 
$$\tilde\theta_{2^,\fE}= \tfrac13(\alpha_1+\alpha_2-2\alpha_3)<0.$$ 
Then $\tilde\theta_{2^,\fE}$ is not nef, so that   $\theta_{2^,\fE}$ is not nef either, 
and the class $\lambda_{2^,\fE}$ is in turn not nef, 
for the argument contained in the proof of Theorem \ref{lnis}.

\subsection{Conclusion of the proof of Theorem \ref{secthm}.} We need to show that if all classes $\lambda_{s,\fE}\in N^1(\PP Q_{s,\fE})$ are nef, then $\fE $ is 
semistable. 

Let   $\fE'=(E',\phi)$ be a rank   $s$ locally-free  Higgs quotient of 
$\fE $. Then   there is a  
section $\sigma\colon C\to 
\hgrass_s(\fE)$  such that $E'=\sigma^\ast Q_s$. Consider the  
curve $C_\sigma=\sigma(C)\subset 
\hgrass_s(\fE)$, the restriction $Q_s^\sigma={Q_s}_{\vert  
C_\sigma}$ and the class 
$\bar\lambda_s={\lambda_{s,\fE}}_{\vert \PP Q_s^\sigma}$. Since  
$\lambda_{s,\fE}$ is nef 
by hypothesis, the class $\bar\lambda_s$ is nef as well. On the other  
hand, we have 
$$ \bar\lambda_s = \lambda_{\PP Q_s^\sigma} + (\mu(E')-\mu(E)) F_s^\sigma$$ 
where $F_s^\sigma$ is the class of the fibre of the projection $\PP  
Q_s^\sigma\to C_\sigma$, 
and $$\lambda_{\PP Q_s^\sigma} = [c_1(\Oc_{\PP  
Q_s^\sigma}(1))]-\mu(Q_s^\sigma)\,F_s^\sigma$$ 
(note that $\mu(Q_s^\sigma)=\mu(E'))$. Since $(\lambda_{\PP  Q_s^\sigma})^s=0$ and $(\lambda_{\PP Q_s^\sigma})^{s-1} 
\cdot F_s^\sigma=1$, the condition $(\bar\lambda_s)^s\ge 0$ implies 
$\mu(E')\ge\mu(E)$, so that $\fE $ is semistable. 

\begin{example} \label{exMiyaoka} \rm Let $X$ be a smooth projective surface 
over $\C$ with ample canonical class $K$. As an application of the criterion established in Theorem \ref{secthm}
we prove the semistability of the Higgs bundle
$ F= \Omega_X\oplus\Oc_X$, with  Higgs structure   given by the morphism
$\phi$ (cf.~\cite{Simp3})
\begin{eqnarray*}\Omega_X\oplus\Oc_X &\to& (\Omega_X\otimes\Omega_X)
\oplus\Omega_X \\
(\omega,f) &\mapsto& (0,\omega)\,.\end{eqnarray*} 
The interest of this example is that since $(F,\phi)$ is semistable
it satisfies the Bogomolov inequality (which holds true also for semistable Higgs bundles, cf. \cite{Simp1}), which in this case yields
the Miyaoka-Yau inequality $3c_2(X)\ge c_1(X)^2$.

Since $K$ is ample, $X$ admits
a K\"ahler-Einstein metric \cite{Y}, hence the cotangent bundle
$\Omega_X$ is semistable with respect to the polarization $K$.
Let $C$ be a curve in the linear system $\vert mK\vert $, with $m$ big enough
for $\Omega=\Omega_{X\vert C}$ to be semistable. Let 
$E= F_{\vert C}= \Omega\oplus\Oc_C$. It is sufficient to prove that
the Higgs bundle $\fE $ is semistable.

By analyzing the possible rank-1 locally-free Higgs quotients of $E$
one finds that $\hgrass_1(\fE)$ has two components, one isomorphic
to $C$, with $Q_1\simeq\Omega_C$ and $\lambda_{1,\fE}=c_1(\Omega_C)$,
which is nef; the other component is isomorphic to $\PP\Omega$ with
$Q_1=\Oc_{\PP\Omega}(1)$ and  $\lambda_{1,\fE}=\lambda_{\PP\Omega}$, which is nef because $\Omega$
is semistable.

For rank-2 locally-free Higgs quotients we find that $\hgrass_2(\fE)$
has two components both isomorphic to $C$. In one case
$Q_2\simeq\Omega$ and $\lambda_{2,\fE}=c_1(\Omega)$,
and in the other $Q_2\simeq\Omega_C\oplus\Oc_C$ with
$\lambda_{2,\fE}=K_C$, which is nef. So $\fE $ is semistable.
\end{example}

\medskip

\section{The higher-dimensional case\label{4}}
In this section we prove Theorem \ref{thirdthm}.
In extending the semistability criterion to the higher-dimensional case
we shall use transcendental techniques. So we assume that $X$ 
is a projective $n$-dimensional smooth variety over $\C$ with
a choice of a polarization $H$. Let 
$\fE=(E,\phi)$  be a rank $r$ Higgs bundle bundle on $X$, and let $\Delta(E)$ be 
the characteristic class
$$\Delta(E) =  c_2(E)-\frac{r-1}{2r}c_1(E)^2 = \frac{1}{2r}c_2(E\otimes E^\ast)\, .$$
We shall denote by $\fE\otimes\fE^\ast$ the Higgs bundle
$(E\otimes E^\ast,\psi)$ where $\psi$ is obtained by coupling
$\phi$ and $\phi^\ast$ in the usual way.

$\hgrass_s(\fE)$   denotes as before the Grassmannian
of Higgs quotients of $\fE$, while $\fQ_s=(Q_s,\Phi_s)$ is the universal
quotient Higgs bundle on $\hgrass_s(\fE)$. Denoting by
$\pi_s\colon \PP Q_s\to X$ the projection onto $X$,
one defines the classes $\lambda_{s,\fE}\in N^1(\PP Q_s)$
as $$\lambda_{s,\fE}=c_1(\OPQ{s}(1))-\frac1r \pi^\ast_s c_1(E).$$

We shall need  a result similar to Proposition \ref{vanishchern} in the case
of Higgs bundles.
We prove a weaker version which is sufficient to our purpose, but very likely
  this result can be strengthened. 
\begin{defin} A Higgs bundle $\fE$ is said to be Higgs-nef
if the line bundle $\OPQ{s}(1)$ is nef for every $s=1,\dots,r-1$ (or, in other terms,
if all bundles $Q_s$ are nef). If both $\fE$ and $\fE^\ast$ are
Higgs-nef, $\fE$ is said to be Higgs-numerically flat.
\end{defin} 
  
\begin{lemma} If $\fE$ is  semistable and Higgs-numerically flat   with $c_1(E)=0$ and there exists a section $\sigma\colon X\hookrightarrow \hgrass_1(\fE)$ of the Higgs Grassmannian $\rho_1\colon \hgrass_1(\fE) \to X$, then all   Chern classes of $E$ vanish.
\label{vanishchernhiggs}\end{lemma}
\begin{proof}  We first prove the statement when $X$ is a surface. By the definition of Higgs-nefness, all the universal bundles $Q_s$ on $\hgrass_s(\fE)$ are nef. Let us consider the exact sequence
$$
0\to S_1 \to \rho_1^\ast E \to Q_1\to 0\,.
$$
Under the identification $\hgrass_{r-1}(\fE^\ast)\simeq \hgrass_1(\fE)$ the bundle $S_1$  gets identified with $Q_{r-1}^\ast$. If we take $\sigma^\ast$ in the above sequence we have
$$
0\to \sigma^\ast(S_1) \to E \to \sigma^\ast(Q_1)\to 0\,.
$$
Now $\sigma^\ast(Q_1)$ and $\sigma^\ast(S_1)^\ast$ are nef. By Theorem 2.5 of \cite{DPS} (but since we are on a surface one can also easily prove this by direct computation) we have $c_1(\sigma^\ast(S_1)^\ast)^2-c_2(\sigma^\ast(S_1)^\ast)\ge 0$. Since $c_1(E)=0$ one has $c_2(E)\le 0$, which together with the Bogomolov inequality $c_2(E)\ge 0$  yields  $c_2(E)=0$ as desired.

If $n=\dim X>2$, taking $m\gg0$ and a smooth hypersurface $Y$ in the linear series $|mH|$, the restriction of $\fE$ to $Y$ is still Higgs-semistable and Higgs-numerically flat with vanishing first Chern class, and   $\hgrass_1(\fE_{|Y})$ has a section. We can iterate this until we get a surface $Z$. Then we have $c_2(E_{|Z})=0$ and therefore  $c_2(E)\cdot  H^{n-2}=0$. Then by
\cite[Thm. 2] {Simp1} 
 $\fE$ is an extension of stable Higgs bundles with vanishing Chern classes, so that the Chern classes of $E$ vanish as well.
\end{proof} 

\subsection*{Proof of Theorem \ref{thirdthm}.} We first prove that i) implies ii), dividing the proof into steps.

Step 1. Let $f\colon C\to X$ be a morphism, where $C$
is a smooth projective curve. Then $f^\ast \fE$ is semistable.
Indeed the definitions of the   bundles $Q_s$ and of the classes
$\lambda_{s,\fE}$ are functorial, so that   
$\lambda_{s,f^\ast\fE}=\tilde f^\ast\lambda_{s,\fE}$, where
$\tilde f\colon \PP Q'_s \to \PP Q_s$ is the map  induced
by $f$ (here $Q'_s$ is the universal quotient Higgs bundle of $f^\ast E$). This implies that the classes $\lambda_{s,f^\ast\fE}$ are
nef.  By Theorem \ref{secthm}  $f^\ast \fE$  is semistable.

Step 2. We prove that $\fE$ is semistable. Indeed by first restricting to
the generic divisor in $\vert mH\vert $ for $m$ big enough, and then
iterating, we may assume that $X$ is a surface. Applying the previous
step to a generic curve in $\vert mH\vert $, again for $m$ big enough, we obtain that
$\fE_{\vert mH}$ is semistable, and then $\fE$ is semistable.

Step 3. We prove that $\fF=\fE\otimes\fE^\ast$ is Higgs-numerically flat. Since $\fF$ is isomorphic to its dual, the point is   to prove that it is Higgs-nef. Fix $s$ and let $\bar f\colon C \to \PP Q_s(\fF)$ be a finite morphism, $C$ being a connected smooth curve.  By Step 1, if $f=\pi_s\circ \bar f$, then $f^\ast\fE$ is semistable, and hence $f^\ast \fF$ is semistable so that $\lambda_{s,f^\ast\fE}$ is a nef class. Since  $c_1(F)=c_1(E\otimes E^\ast)=0$, the class  $\lambda_{s,f^\ast\fE}$ equals $c_1(\cO_{\PP f^\ast(Q_s(\fF))}(1))$.  Let $\tau\colon \bar f(C) \to \PP Q_s(\fF)$ be the section induced by $\bar f$; one has
$$
\deg \tau^\ast \cO_{\PP f^\ast(Q_s(\fF))} (1)= [\bar f(C)]\cdot 
c_1(\cO_{\PP (Q_s(\fF))}(1))
$$
so that $Q_s(\fF)$ is nef, i.e., $\fE$ is Higgs-nef.

Step 4. We  show that $\Delta(E)=0$. Since $\fF$ is Higgs-numerically flat and semistable as a Higgs bundle and the Higgs Grassmannian $\hgrass_1(\fF)$ has a section induced by the evaluation morphism $\fF\to \cO_X$, all   Chern classes of $E\otimes E^\ast$ vanish by Lemma \ref{vanishchernhiggs}, whence the claim.

Now we prove the converse statement, again dividing it into steps.

Step 1. The present hypotheses imply
that $\fF$ is semistable, and $c_2(E\otimes E^\ast)=0$, while
of course $c_1(E\otimes E^\ast)=0$. By \cite[Thm. 2]{Simp1}  there is   a filtration in  Higgs subbundles
$$ 0 = \fF_0 \subset \fF_1\subset\dots\subset \fF_s=\fF$$
such that every quotient $\fF_i/\fF_{i-1}$ is stable and has vanishing Chern classes.
Again by results contained in  \cite{Simp1} we know that
each quotient $\fG_i=\fF_i/\fF_{i-1}$ admits a Hermitian-Yang-Mills metric.
Let $\Omega_i$ be the curvature of the associated Chern connection.
Since $c_1(G_i)=c_2(G_i)=0$, we have
$$ 0 = \int_X \mbox{tr}(\Omega_i\wedge \Omega_i) \cdot H^{n-2} =
\gamma_1 \Vert \Omega_i\Vert ^2 - \gamma_2\Vert\Lambda \Omega_i\Vert^2 = \gamma_1 \Vert \Omega_i\Vert ^2$$
for some positive constants $\gamma_1$, $\gamma_2$, so that the Chern connection of $G_i$ is flat.

Step 2. For a fixed $s$ with $0<s<r$, let
$\bar f\colon C\to \PP Q_s$ be a finite morphism, where 
 $C$ is a smooth irreducible projective curve. 
Let $f\colon C \to X$ be the composition
of $\bar f$ with the projection $\PP Q_s\to X$. We show that the Higgs bundle $f ^\ast \fE$ is semistable. Indeed the Higgs bundle  $f ^\ast\fF$ is filtered by the  Higgs bundles $f ^\ast \fF_i$, and the pullbacks $f ^\ast(\fF_i/\fF_{i-1})\simeq f ^\ast \fF_i/f ^\ast \fF_{i-1}$ carry flat unitary connections, hence they are polystable (again \cite[Thm. 1]{Simp1}). Moreover they all have degree zero. As a consequence, $f ^\ast \fF$ is semistable, and $f ^\ast \fE$ is semistable as well.
 
 Step 3. By  Theorem \ref{secthm}, the classes $\tilde\lambda_{s,\fE}$ in
 $N^1(\PP Q'_s)$ are  nef. This implies that for any irreducible
 curve $C'\subset \PP Q_s$ one has  $[C'] \cdot \lambda_{s,\fE} \ge 0$,
 i.e., it implies that  all classes $\lambda_{s,\fE}$ are nef. Indeed if $C$ is the normalization of $C'$ we can apply the previous constructions to $C$.

This concludes the proof of Theorem \ref{thirdthm}.

Theorem  \ref{thirdthm} has an immediate Corollary.

\begin{corol} A semistable Higgs bundle $\fE=(E,\phi)$ on an $n$-dimensional
projective polarized complex manifold $(X,H)$ such that
$c_1(E)\cdot H^{n-1}=\mbox{\rm ch}_2(E)\cdot H^{n-2} = 0$ is Higgs-numerically
flat. \end{corol}
\begin{proof} Again by \cite[Thm. 2]{Simp1}  all Chern classes of $E$
vanish. So $\Delta(E)=0$, and by  Theorem  \ref{thirdthm}, all classes
$\lambda_{s,\fE}$ are nef. But since $c_1(E)=0$ this implies that
all bundles $Q_s$ are nef, i.e., $\fE$ is Higgs-nef. Applying the same argument
to the dual Higgs bundle $\fE^\ast$ one obtains the claim.
\end{proof}

\end{document}